\documentstyle[11pt]{article}
\newcommand{\ineqS}{$\mu_S(I)\mu_S(S-I)>\mu_S(I+(S-I))$}

\newcommand{\Z}{{\bf Z}}

\begin{document}

\vspace{0.5in}

\title{Perfect Pairs of Ideals and Duals in Numerical Semigroups}  
\bigskip

\author{Kurt Herzinger\\
\ \\
kurt.herzinger@usafa.af.mil\\
\ \\
\newline
Department of Mathematical Sciences\\
2354 Fairchild Dr Suite 6D124\\  
United States Air Force Academy\\
USAF Academy, CO 80840-6252\\
Fax: 719-333-2114\\
\ \\
\ \\
Stephen Wilson \hspace{1.2in} N\'andor Sieben \\
\ \\
Stephen.Wilson@nau.edu \ \ \ \ \ Nandor.Sieben@nau.edu\\
\ \\
Jeff Rushall \\
\ \\
Jeffrey.Rushall@nau.edu\\
\newline
\ \\
Department of Mathematics and Statistics\\
Northern Arizona University\\
Flagstaff, AZ 86011-5717\\
Fax: 928-523-5847}
\date{}
\maketitle



\parbox{4.0in}
{{\bf Abstract.} \em This paper considers numerical semigroups $S$ that have a non-principal relative 
ideal $I$ such that $\mu_S(I)\mu_S(S-I)=\mu_S(I+(S-I)) $.  We show the existence of an infinite family of such pairs $(S,I)$ in 
which $I+(S-I)=S\backslash\{0\}$.  We also show examples of such pairs that are not members of this 
family.  We discuss the computational process used to find these examples and present some open questions 
pertaining to them.}
\vspace{.2in}

%%%%%%%%%%%%%%%%%%   Section 1   %%%%%%%%%%%%%%%%%%

\begin{center}
{\bf 1. Definitions, Notation and Background}\\
\end{center}

\noindent
{\bf (1.1) Definitions/Notations:}\\

\noindent
(a) A {\it numerical semigroup} $S$ is a subset of the non-negative integers {\bf N} which contains $0$, is closed under addition, and such that {\bf N}$\backslash S$ is finite.  If $G$ is the smallest subset of $S$ such that every element of $S$ is a sum of elements from $G$, then we say $G$ is {\it the minimal generating set of} $S$ and we write $S = \langle G \rangle$.\\

\noindent
(b) The {\it multiplicity} of $S$, denoted by $m(S)$, is the smallest positive element of $S$.\\

\noindent
(c) $g(S) =$ max$({\bf N} \backslash S)$ is called the {\it Frobenius number} of $S$.\\

\noindent
(d) $n(S)$ represents the number of elements of $S$ that are less than $g(S)$.\\

\noindent
(e) We say that $S$ is {\it symmetric} provided $g(S)$ is odd and $n(S)=\frac{g(S)+1}{2}$.  Several equivalent definitions can be found in \cite{BDF}.\\

\noindent
(f) A {\it relative ideal} $I$ of $S$ is a finite union of cosets $z+S$ where $z \in {\bf Z}$.  The notation $I = (z_1,\dots , z_k)$ means $I = (z_1+S) \cup \dots \cup (z_k+S)$ and $I$ cannot be written as a union of cosets by any proper subset of $\{z_1,\dots , z_k\}$.  We refer to $\{z_1,\dots , z_k\}$ as {\it the minimal generating set of} $I$.  \\

\noindent
(g) $\mu_S(I)$ represents the size of the minimal generating set for $I$.\\

\noindent
(h) The {\it dual} of $I$ in $S$ is $S-I=\{z \in \Z\ \mid z+I\subseteq S\}$.\\

\noindent
(i) If $I$ and $J$ are relative ideals of $S$, we define their sum by $I+J =\{i+j \mid i \in I , j \in J\}$.\\
  
\noindent
(j) We define the {\it Apery set of} $S$ {\it with respect to} $m(S)$ by $Ap(S) = \{s \in S \mid s-m(S) \notin S\}$.  The properties of $Ap(S)$ and some of its subsets are detailed in \cite{MCH}.\\

\vspace{.1in}

\noindent
{\bf (1.2) Example:}  Let $S = \langle 10,11,13,17,19\rangle$ and $I = (2,5)$.  That is, $S$ has minimal generating set $\{10,11,13,17,19\}$ and $I = (2+S) \cup (5+S)$.  Listing the elements of $S$ we find

\begin{center}

$S = \{0,10,11,13,17,19,20,21,22,23,24,26 \rightarrow \}$\\

\end{center}

\noindent
where `$\rightarrow$' indicates that all integers greater than 26 are in $S$.  We see that $m(S) = 10$, $g(S) = 25$, and $n(S) = 11$ so $S$ is not symmetric.  Further, we have $Ap(S) = \{0,11,13,17,19,22,24,26,28,35\}$.\\

Listing the elements of $I$ we have

\begin{center}

$I = \{2,5,12,13,15,16,18,19,21,\rightarrow \}$\\

\end{center}

\noindent
Computing $S-I$ we find

\begin{center}

$S-I = \{8,15,17,18,19,21,22,24,\rightarrow \}$\\

\end{center}

\noindent
Extracting the minimal generating set for $S-I$ yields

\begin{center}

$S-I = (8,15,17,22,24)$\\

\end{center}

\noindent
At this point we have $\mu_S(I) = 2$ and $\mu_S(S-I)=5$.  We now form a generating set for $I+(S-I)$ 
by adding each of the generators of $I$ to each of the generators of $S-I$.  We find $\{10,13,17,20,19,22,24,27,26,29\}$ is 
a generating set for $I+(S-I)$.  A quick check reveals that this generating set is not minimal 
because $29-10 \in S$, $26-13 \in S$ etc.  We discover

\begin{center}

$I+(S-I)=(10,13,17,19,22)$\\

\end{center}

\noindent
so $\mu_S(I+(S-I))=5$ and hence \ineqS.\\

\pagebreak

\noindent
{\bf Background and Motivation:} Let $S$ be a numerical semigroup and let $I$ be a non-principal relative ideal of $S$.  The inequality 
$$\mu_S(I)\mu_S(S-I)\ge\mu_S(I+(S-I))\ \ \ (*)$$ 
always holds and several conditions that imply strict inequality have been established.  For example, if $m(S) \le 8$, then \ineqS\ always holds (see \cite{HS}).  Further, if we restrict to the case $\mu_S(I)=\mu_S(S-I)=2$, then $\mu_S(I+(S-I)) \le 3$ whenever $m(S) \le 9$ (see \cite{H3}).  Until now, the investigation of the inequality $(*)$ has revealed only three examples of numerical semigroups $S$ and non-principal relative ideals $I$ such that equality holds (see \cite{HS}).

In this first section we introduce the notions of {\it balanced} and {\it unitary} numerical semigroups and establish some basic properties related to these structures. In section 2 we prove that unitary numerical semigroups $S$ have a non-principal relative ideal $I$ such that $\mu_S(I)=2$, $\mu_S(S-I)=2$, $\mu_S(I+(S-I))=4$ so that the equality in $(*)$ holds, and that $I+(S-I)=S\backslash \{0\}$.  We further prove that each $S$ in this family is symmetric and establish the value of the Frobenius number for each $S$.  In section 3 we discuss how the search for these examples was conducted.  We will also provide examples of the equality in $(*)$ that are different from those described above.  We conclude with some open questions related to this investigation.

Suggested background reading for numerical semigroups and their connections to commutative algebra include \cite{BDF}, \cite{FGH} and \cite{K}.  
The original motivation for the investigation in this paper comes from the study of 
torsion in tensor products of modules over certain types of rings.  The specifics of the relationship between this topic and numerical semigroups are detailed in \cite{H1}. Details concerning the investigation of torsion in tensor products can be found in \cite{A}, \cite{C}, \cite{HW1} and \cite{HW2}.\\

\vspace{.07in}

\noindent
{\bf (1.3) Definitions/Notation:}  Let $B = \{a_1,a_2,a_3,a_4\}$ be a set of four positive integers satisfying the following conditions:

\begin{center}
\begin{tabular}{rl}
(1) & $a_1<a_2<a_3<a_4$ \\
(2) & ${\rm gcd}(a_1,a_2,a_3,a_4)=1$\\
(3) & $a_i \not|\ a_j$ for $i \not= j$\\
(4) & $a_1+a_4 = a_2+a_3$\\

\end{tabular}
\end{center}

We call such a set {\it balanced}.  We refer to the quantity in (4) as the {\it common sum} of $B$, denoted by $CS(B)$.  Note that $CS(B)$ is divisible by both ${\rm gcd}(a_1,a_4)$ and ${\rm gcd}(a_2,a_3)$.  Since these gcd's are relatively prime by (2), we know that $CS(B)$ is divisible by their product.  We define the {\it common quotient} of $B$ to be the integer $$CQ(B) = \frac{CS(B)}{{\rm gcd}(a_1,a_4)\cdot {\rm gcd}(a_2,a_3)}.$$  We say a balanced set $B$ is {\it unitary} provided $CQ(B) = 1$.  Finally, we say a numerical semigroup $S$ is {\it balanced} ({\it unitary}) provided $S$ is minimally generated by a balanced (unitary) set.  In this case the symbols $CS(S)$ and $CQ(S)$ refer to the common sum and the common quotient of the generating set of $S$. \\

\noindent
{\bf (1.4) Examples:}\\

\noindent
(1) $S = \langle 14,15,20,21\rangle$ is an example of a unitary numerical semigroup.\\

\noindent
(2) $S = \langle12,15,25,28\rangle$ is an example of a balanced numerical semigroup that is not unitary since $CQ(S) = 2$.\\

\noindent
(3) $S = \langle10,14,15,21\rangle$ is an example of a numerical semigroup that is not balanced.\\

Note that the terms {\it balanced} and {\it unitary} are reserved for sets of exactly four elements.\\

\noindent
{\bf (1.5) Notation:}\ Let $S = \langle a_1,a_2,a_3,a_4\rangle$ be a balanced numerical semigroup.  Let $D = {\rm gcd}(a_1,a_4)$ and $E = {\rm gcd}(a_2,a_3)$.  We will then write $$a_1 = q_1D,\ a_2 = q_2E,\ a_3 = q_3E,\ a_4 = q_4D.$$\\

\noindent
{\bf (1.6) Examples:}\\

\noindent
(1) If $S = \langle14,15,20,21\rangle$, then $q_1 = 2$, $q_2 = 3$, $q_3 = 4$, $q_4 = 3$, $D = 7$, and $E = 5$.\\

\noindent
(2) If $S = \langle12,15,25,28\rangle$, then $q_1 = 3$, $q_2 = 3$, $q_3 = 5$, $q_4 = 7$, $D = 4$, and $E = 5$.\\

\pagebreak

\noindent
{\bf (1.7) Notes:}\\

\noindent
(1) ${\rm gcd}(q_1,q_4)={\rm gcd}(q_2,q_3)={\rm gcd}(D,E)={\rm gcd}(q_1,E)={\rm gcd}(q_2,D)={\rm gcd}(q_3,D)={\rm gcd}(q_4,E)=1$\\

\noindent
(2) The assumption $a_1<a_2<a_3<a_4$ implies $q_1<q_4$ and $q_2<q_3$.\\

\noindent
(3) If $q_1=1$ or $q_2=1$, then $a_1\mid a_4$ or $a_2\mid a_3$.  Therefore, $q_1>1$ and $q_2>1$.\\

\noindent
(4) $S$ is unitary if and only if $q_2+q_3 = D$ if and only if $q_1+q_4 = E$ if and only if $CS(S)=DE$.\\

\noindent
(5) ${\rm lcm}(a_2,a_3)=q_3a_2 = q_2a_3$ and ${\rm lcm}(a_1,a_4)=q_4a_1 = q_1a_4$.\\

\vspace{.2in}

At this point we examine the Apery set of a balanced numerical semigroup.  Recall that $Ap(S) = \{s \in S \mid s - m(S) \notin S\}$.  If $S=\langle a_1,a_2,a_3,a_4\rangle$, then every element $s \in S$ can be expressed in the form $s = t_1a_1+t_2a_2+t_3a_3+t_4a_4$ where $t_i \ge 0$ for each $i$.  Now assume that $S$ is balanced and $s \in Ap(S)$.  Consider the following: \\

($i$) If $t_1 \ge 1$, then $s-a_1 = (t_1-1)a_1+t_2a_2+t_3a_3+t_4a_4 \in S$.  Since $a_1 = m(S)$, we see that $s \notin Ap(S)$, a contradiction.\\

($ii$) If $t_4 \ge q_1$, then by (1.7 (5)) we have $s = t_2a_2+t_3a_3+(t_4-q_1)a_4+q_1a_4 = q_4a_1+t_2a_2+t_3a_3+(t_4-q_1)a_4 \notin Ap(S)$ by ($i$), a contradiction.\\

($iii$) If $t_2 > 0$ and $t_3 > 0$, then by (1.3 (4))  $s = (t_2-1)a_2+(t_3-1)a_3+t_4a_4+a_2+a_3 = a_1+(t_2-1)a_2+(t_3-1)a_3+(t_4+1)a_4 \notin Ap(S)$ by ($i$), a contradiction.\\

($iv$) If $t_2>q_3$, then we can use (1.7 (5)) to write $s = (t_2-q_3)a_2+q_3a_2+t_4a_4 = (t_2-q_3)a_2+q_2a_3+t_4a_4$.  We conclude $s \notin Ap(S)$ by ($iii$), a contradiction.  Similarly if $t_3>q_2$, then $s \notin Ap(S)$.\\

\vspace{.1in}

This discussion is summarized in the following proposition.\\

\pagebreak

\noindent
{\bf (1.8) Proposition:}\ If $S$ is a balanced numerical semigroup and $s \in Ap(S)$, then 

\begin{tabular}{rcl}
$s = t_4a_4$ & where & $0 \le t_4 \le q_1-1$\ \ \ \ or \\
$s = t_2a_2+t_4a_4$ & where & $1 \le t_2 \le q_3$ and $0 \le t_4 \le q_1-1$\ \ \ \ or\\
$s = t_3a_3+t_4a_4$ & where & $1 \le t_3 \le q_2-1$ and $0 \le t_4 \le q_1-1$.\\

\end{tabular}

\vspace{.2in}

The following proposition precisely identifies the elements in the Apery set of a unitary numerical semigroup.\\

\noindent
{\bf (1.9) Proposition:}\ Let $S$ be a balanced numerical semigroup.  Define the following sets:

\hspace{.8in}$A_1 = \{t_4a_4 \mid 0 \le t_4 \le q_1-1\}$ 

\hspace{.8in}$A_2 = \{t_2a_2+t_4a_4 \mid 1 \le t_2 \le q_3 , 0 \le t_4 \le q_1-1\}$

\hspace{.8in}$A_3 = \{t_3a_3+t_4a_4 \mid 1 \le t_3 \le q_2-1 , 0 \le t_4 \le q_1-1\}$.\\

If $S$ is unitary, then  $Ap(S) = A_1 \cup A_2 \cup A_3$.  Moreover, $A_1$, $A_2$ and $A_3$ are pairwise disjoint and the elements of $A_1 \cup A_2 \cup A_3$ are pairwise non-congruent modulo $a_1$.\\

\noindent
{\bf Proof:}\ From basic set theory we know $|A_1 \cup A_2 \cup A_3| \le |A_1|+|A_2|+|A_3| \le q_1 + q_3q_1 + (q_2-1)q_1 = q_1(1+q_3+q_2-1) = q_1(q_3+q_2) = q_1D\cdot CQ(S)$.  If $S$ is unitary, we then have $|A_1 \cup A_2 \cup A_3| \le q_1D = a_1$.  Now (1.8) says $Ap(S) \subseteq A_1 \cup A_2 \cup A_3$.  Since $|Ap(S)| = m(S) = a_1$, we conclude $Ap(S) = A_1 \cup A_2 \cup A_3$.  The other statements of the proposition now follow from the fact that the elements of $Ap(S)$ are pairwise non-congruent modulo $a_1$.\\

As the following example shows, if $S$ is balanced but not unitary, the result of (1.9) may not hold.\\

\noindent
{\bf (1.10) Example:}\ If $S = \langle12,15,25,28\rangle$, then $2a_3+2a_4 = 106 \in A_3$ but is not an element of $Ap(S)$ because $106 - a_1 = 94 \in S$.\\

We finish this section by introducing the notion of a $k \times m$ brick.\\

\noindent
{\bf (1.11) Definition:}\ Let $S$ be a numerical semigroup and let $I$ be a relative ideal of $S$.  We refer to the pair $(S,I)$ as a $k \times m$ {\it brick} provided $\mu_S(I) = k$ and $\mu_S(S-I)=m$ and $\mu_S(I+(S-I)) = \mu_S(I)\mu_S(S-I)=km$.  Further, if $I+(S-I) = S\backslash \{0\}$ we say that $(S,I)$ is a {\it perfect} $k \times m$ {\it brick}.\\

\vspace{.2in}

At the end of \cite{HS} it is stated that there are only three known examples of bricks:\\

\begin{center}
\begin{tabular}{cl}

(1) & $(\langle 14,15,20,21\rangle,(0,1))$ \\
(2) & $(\langle 10,14,15,21\rangle,(0,1))$ \\
(3) & $(\langle 14,15,20,21,25\rangle,(0,1))$ \\

\end{tabular}
\end{center}

\vspace{.2in}

All three examples are $2 \times 2$ bricks and (1) is perfect while (2) and (3) are not.  In the next section we will prove the existence of infinitely many perfect $2 \times 2$ bricks.\\

%%%%%%%%%%%%%%%%%%   Section 2   %%%%%%%%%%%%%%%%%%

\begin{center}

{\bf 2. Unitary Numerical Semigroups Yield Perfect Bricks}\\

\end{center}

We state the main result of this investigation in the following theorem.\\

\noindent
{\bf (2.1) Theorem:}\ Let $S = \langle a_1,a_2,a_3,a_4\rangle$ be a unitary numerical semigroup and let $n = a_2-a_1=a_4-a_3$.  If $I=(0,n)$, then $(S,I)$ is a perfect $2 \times 2$ brick.\\

The proof of this theorem will be established through a series of lemmas and propositions.  In (2.6) we prove that unitary numerical semigroups are symmetric and a formula for the Frobenius number will be given.  In (2.8) we show that the dual of the relative ideal $I$ in the statement of the theorem is the relative ideal $(a_1,a_3)$.  It follows that $I+(S-I) = (a_1,a_2,a_3,a_4)$ and hence $(S,I)$ is a perfect $2 \times 2$ brick.  We will continue to use the assumptions and notations from section 1.\\

Let $S$ be a unitary numerical semigroup.  We begin by examining a subsemigroup of $S$.  Let $T = \langle a_1,a_2,a_3\rangle$.  We will prove that $T$ is symmetric and establish a formula for $g(T)$.  This will lead us to analogous results for $S$.\\

Let $d_1 = {\rm gcd}(a_2,a_3) = E$, $d_2 = {\rm gcd}(a_1,a_3) = {\rm gcd}(q_1,q_3)$ and $d_3 = {\rm gcd}(a_1,a_2) = {\rm gcd}(q_1,q_2)$.  The derived semigroup of $T$ (see \cite{FGH}) is\\

\begin{center}
\begin{tabular}{rl}

$\overline{T}$\hspace{-.08in} &$=\langle\frac{a_1}{d_2d_3},\frac{a_2}{d_1d_3},\frac{a_3}{d_1d_2}\rangle$ \vspace{.06in}\\
  &$=\langle\frac{q_1D}{d_2d_3},\frac{q_2}{d_3},\frac{q_3}{d_2}\rangle$.\\
  
\end{tabular}
\end{center}
\ \\
Note that $d_2\mid q_1$ and $d_3\mid q_1$ and ${\rm gcd}(d_2,d_3)=1$ so $d_2d_3\mid q_1$.  Define $t$ by $d_2d_3t=q_1$.  Now,\\

\begin{center}
\begin{tabular}{rll}

$d_3t(\frac{q_2}{d_3}) + d_2t(\frac{q_3}{d_2})$\hspace{-.08in} &$=\frac{d_2d_3tq_2+d_2d_3tq_3}{d_2d_3}$ & \vspace{.06in}\\
  &$=\frac{d_2d_3t(q_2+q_3)}{d_2d_3}$ & \vspace{.06in}\\
  &$=\frac{q_1D}{d_2d_3}$. & ($\dagger$)\\
  
\end{tabular}
\end{center}

We conclude that $\frac{q_1D}{d_2d_3} \in \langle\frac{q_2}{d_3},\frac{q_3}{d_2}\rangle$ whence $\overline{T}$ is minimally generated by 2 elements.  Therefore $T$ is symmetric by \cite{FGH} (corollary to Theorem 14).\\

Next, from ($\dagger$) we have $d_2d_3tq_2+d_2d_3tq_3 = a_1$.  This says $$a_1 \in \langle q_2,q_3 \rangle = \langle \frac{a_2}{d_1},\frac{a_3}{d_1} \rangle. \ \ (\dagger\dagger)$$ 

From \cite{S}, we know that for relatively prime positive integers $\alpha$ and $\beta$ we have $$g(\langle \alpha,\beta \rangle) = \alpha\beta-\alpha-\beta . \ \ (\dagger\dagger\dagger)$$\\

\vspace{-.3in}
We can now compute the value of $g(T)$.  The first equality in the following manipulation comes from (\cite{FGH}, Prop. 8(c)) and the fact that ${\rm gcd}(a_1,d_1)={\rm gcd}(a_1,E)=1$ (see also \cite{J}). \\

\begin{tabular}{lll}

$g(T)$\hspace{-.08in} &=\ $d_1g(\langle a_1,\frac{a_2}{d_1},\frac{a_3}{d_1}\rangle)+(d_1-1)a_1$ & \\
\ \ \ \ \ &=\ $d_1g(\langle\frac{a_2}{d_1},\frac{a_3}{d_1}\rangle)+(d_1-1)a_1$ & (by $\dagger\dagger$)\\
 &=\ $d_1(\frac{a_2}{d_1}\frac{a_3}{d_1}-\frac{a_2}{d_1}-\frac{a_3}{d_1})+(d_1-1)a_1$ & (by $\dagger\dagger\dagger$)\\
 &=\ $q_2a_3 + Ea_1 - a_1 - a_2 - a_3$ & \\
 &=\ $(q_2-1)a_3 + (E-1)a_1 - a_2$ & \\
 &=\ $(q_2-1)a_3 + (E-1-q_4)a_1 + q_4a_1 - a_2$ & \\
 &=\ $(q_2-1)a_3 + (q_1-1)a_1 + q_4a_1 - a_2$ & (by (1.7 (4)) \\
 &=\ $(q_2-1)a_3 + (q_1-1)a_1 + q_1a_4 - a_2$ & (by (1.7 (5)) \\
 &=\ $(q_2-1)a_3 + (q_1-1)a_1 + q_1a_4 - a_1 - a_4 + a_3$ & (by (1.3 (4)) \\
 &=\ $q_2a_3 + (q_1-2)a_1 + (q_1-1)a_4$ & \\
 &=\ $q_3a_2 + (q_1-2)a_1 + (q_1-1)a_4$ & (by (1.7 (5)).\\
  
\end{tabular}

\ \\

We summarize the above analysis in the following lemma.\\

\noindent
{\bf (2.2) Lemma:}\ Let $S = \langle a_1,a_2,a_3,a_4\rangle$ be a unitary numerical semigroup and let $T = \langle a_1,a_2,a_3\rangle$.  Then $T$ is symmetric and $g(T) = (q_1-2)a_1 + q_2a_3 + (q_1-1)a_4 = (q_1-2)a_1 + q_3a_2 + (q_1-1)a_4$.\\

\vspace{.1in}

\noindent
{\bf (2.3) Example:}\ Recall from (1.6 (1)), if $S = \langle 14,15,20,21\rangle$, then $q_1 = 2$, $q_2 = 3$ and $q_3 = 4$.  It is quick to confirm that $T=\langle 14,15,20\rangle$ is symmetric and $g(T) = (q_1-2)a_1 + q_2a_3 + (q_1-1)a_4 = 0 + 60 + 21 = 81$.\\

\vspace{.1in}

\noindent
{\bf (2.4) Lemma:}\ Let $S = \langle a_1,a_2,a_3,a_4\rangle$ be a unitary numerical semigroup and let $T = \langle a_1,a_2,a_3\rangle$.  Then $g(S) \ge g(T)-(q_1-1)a_1$.\\

\noindent
{\bf Proof:}\ It suffices to show $g(T)-(q_1-1)a_1 \notin S$.  By way of contradiction, assume $g(T)-(q_1-1)a_1 \in S$.  Then we can write $g(T)-(q_1-1)a_1 = t_1a_1 + t_2a_2 + t_3a_3 + t_4a_4$ where each $t_i \ge 0$.  By the division algorithm we can write $t_4 = zq_1 + r$ where $z \ge 0$ and $0 \le r < q_1$.  Now,\\

\begin{tabular}{rll}

\hspace{-.22in}$g(T)$\hspace{-.08in} &=\ $(q_1-1)a_1+t_1a_1 + t_2a_2 + t_3a_3 + t_4a_4$ & \\
                  &=\ $((q_1-1)+t_1)a_1 + t_2a_2 + t_3a_3 + (zq_1+r)a_4$ & \\
                  &=\ $((q_1-1)+t_1)a_1 + t_2a_2 + t_3a_3 + zq_1a_4 + ra_4$ & \\
                  &=\ $((q_1-1)+t_1)a_1 + t_2a_2 + t_3a_3 + zq_4a_1 + ra_4$ & (1.7 (5))\\
                  &=\ $((q_1-1)+t_1+zq_4)a_1 + t_2a_2 + t_3a_3 + ra_4$ & \\
                  &=\ $((q_1-1)+t_1+zq_4)a_1 + t_2a_2 + t_3a_3 + ra_2 + ra_3 - ra_1$ & (1.3 (4))\\
                  &=\ $((q_1-1)-r+t_1+zq_4)a_1 + (t_2+r)a_2 + (t_3+r)a_3$. & \\
  
\end{tabular}

\vspace{.2in}
The coefficients on $a_1$, $a_2$ and $a_3$ are all non-negative, and so the expression on the right side is an element of $T$.  We conclude $g(T) \in T$, a contradiction.  Therefore, $g(T) - (q_1-1)a_1 \notin S$ and the desired result follows.\\

\vspace{.1in}

\noindent
{\bf (2.5) Lemma:}\ Let $S = \langle a_1,a_2,a_3,a_4\rangle$ be a unitary numerical semigroup and let $T = \langle a_1,a_2,a_3\rangle$.  Then $|S \backslash T| \ge \frac{q_1-1}{2}a_1$.\\

\pagebreak

\noindent
{\bf Proof:}\ Define the following subsets of $S$:\\

\noindent
\begin{tabular}{ll} 
   
$B_1$\hspace{-.1in} & $= \{t_1a_1 + t_4a_4 \mid 1 \le t_4 \le q_1-1 ,0 \le t_1 \le t_4-1\}$ \vspace{.1in} \\

$B_2$\hspace{-.1in} & $= \{t_1a_1+t_2a_2+t_4a_4 \mid 1\le t_2\le q_3 , 1 \le t_4 \le q_1-1 ,0 \le t_1 \le t_4-1\}$ \vspace{.1in} \\

$B_3$\hspace{-.1in} & $= \{t_1a_1+t_3a_3+t_4a_4 \mid 1\le t_3\le q_2-1 , 1 \le t_4 \le q_1-1 , 0 \le t_1 \le t_4-1\}$ \\

\ \\
               
\end{tabular}

First note that if $x \in B_1$, then by (2.2) we have \\

\begin{tabular}{rl}
$g(T)-x$\hspace{-.1in} &=\ $((q_1-2)a_1+q_3a_2+(q_1-1)a_4)-(t_1a_1 + t_4a_4)$\\
         &=\ $(q_1-2-t_1)a_1+q_3a_2+(q_1-1-t_4)a_4$.\\
\end{tabular}

\vspace{.2in}

Since $t_1 \le t_4-1$, we see that $q_1-2-t_1 \ge q_1-1-t_4$.  We can now use (1.3 (4)) to write $g(T) - x = u_1a_1 + u_2a_2 + u_3a_3$ where each $u_i \ge 0$.  We conclude that $g(T)-x \in T$ for all $x \in B_1$.  Similar arguments yield the same result if $x \in B_2$ or $x \in B_3$. Therefore, if $x \in B_1 \cup B_2 \cup B_3$, then $g(T) - x \in T$.  It follows that $x \notin T$ and hence $B_1 \cup B_2 \cup B_3 \subseteq S \backslash T$.  By (1.9) we know $B_1 \cap B_2 = B_1 \cap B_3 = B_2 \cap B_3 = \emptyset$.   Therefore,\\

\begin{tabular}{rl} 
   
$|B_1\cup B_2 \cup B_3|$\hspace{-.08in} & $= |B_1|+|B_2|+|B_3|$ \vspace{.1in}\\
             & $=\frac{(q_1-1)q_1}{2} + \frac{(q_1-1)q_1}{2}q_3 + \frac{(q_1-1)q_1}{2}(q_2-1)$\vspace{.1in}\\
             & $=\frac{(q_1-1)q_1}{2}(1+q_3+q_2-1)$\vspace{.1in}\\
             & $=\frac{q_1-1}{2}q_1D$ \ \ \ \ \ \ \ \ \ \ \ (since $S$ is unitary) \vspace{.1in}\\
             & $=\frac{q_1-1}{2}a_1$\\
              
\end{tabular}

\vspace{.2in}

We conclude $|S \backslash T| \ge \frac{q_1-1}{2}a_1$.\\

\noindent
{\bf (2.6) Proposition:}\ Let $S = \langle a_1,a_2,a_3,a_4\rangle$ be a unitary numerical semigroup and let $T = \langle a_1,a_2,a_3\rangle$.  Then $S$ is symmetric and $g(S) = g(T)-(q_1-1)a_1$.\\

\noindent
{\bf Proof:}\ Define the following sets:
\vspace{.1in}
\begin{center}
\begin{tabular}{rl} 
   
$C_1$\hspace{-.1in} & $= \{t \in T \mid g(S) < t < g(T) \}$ \vspace{.05in}\\
$C_2$\hspace{-.1in} & $= \{x \in S \backslash T \mid x < g(S) \}$ \vspace{.05in} \\
$C_3$\hspace{-.1in} & $= \{y \in S \backslash T \mid y > g(S) \}$ \\
             
\end{tabular}
\end{center}

\pagebreak

Note that $C_1 \cup C_3 = \{g(S)+1,g(S)+2,\dots , g(T)-1,g(T)\}$ and $|C_1|+|C_3|=|C_1 \cup C_3| = g(T)-g(S)$.

  Also note that $C_2 \cup C_3 = S \backslash T$  so $|C_2|+|C_3|=|C_2 \cup C_3|=|S \backslash T| \ge \frac{q_1-1}{2}a_1$ by (2.5).  Now,

\vspace{.08in}

\begin{tabular}{ll} 
   
\hspace{.2in}$\frac{g(S)+1}{2}$\hspace{-.08in} & $\ge n(S)$\vspace{.05in}\\
 & $= n(T) - |C_1|+|C_2|$ \vspace{.05in}\\
 & $= n(T) - |C_1| -|C_3|+|C_2|+|C_3|$\vspace{.05in}\\
 & $= n(T) - (|C_1|+|C_3|)+(|C_2|+|C_3|)$\vspace{.05in}\\
 & $= n(T) - (g(T)-g(S))+(|C_2|+|C_3|)$\vspace{.05in}\\
 & $\ge n(T) - g(T)+g(S)+ \frac{q_1-1}{2}a_1$\vspace{.05in}\\
 & $= \frac{g(T)+1}{2} - g(T)+g(S)+ \frac{q_1-1}{2}a_1$ \ \ \ \ \ ($T$ is symmetric)\vspace{.05in}\\
 & $= \frac{g(T)+1-2g(T)+2g(S)+(q_1-1)a_1}{2}$ \vspace{.05in}\\
 & $= \frac{g(S)+(g(S)-g(T))+1+(q_1-1)a_1}{2}$ \vspace{.05in}\\
 & $\ge \frac{g(S)-(q_1-1)a_1+1+(q_1-1)a_1}{2}$\ \ \ \ \ (by (2.4))\vspace{.05in}\\
 & $= \frac{g(S)+1}{2}$.\\
                
\end{tabular}

\vspace{.2in}

It follows that all of the above expressions are equivalent.  Therefore $g(S) = g(T)-(q_1-1)a_1$ and $n(S)= \frac{g(S)+1}{2}$ which says $S$ is symmetric.\\

The proof of the following statement is clear.\\

\noindent
{\bf (2.7) Corollary:}\ Let $S = \langle a_1,a_2,a_3,a_4\rangle$ be a unitary numerical semigroup and let $T = \langle a_1,a_2,a_3\rangle$.  If $B_1$, $B_2$ and $B_3$ are as defined in the proof of (2.5), then $S \backslash T = B_1 \cup B_2 \cup B_3$.\\

\vspace{.1in}

\noindent
{\bf (2.8) Proposition:}\ Let $S = \langle a_1,a_2,a_3,a_4\rangle$ be a unitary numerical semigroup and let $n = a_2-a_1=a_4-a_3$.  If $I=(0,n)$, then $S-I=(a_1,a_3)$.\\

\noindent
{\bf Proof:}\ It is clear from the definitions of $S-I$ and $n$ that $a_1,a_3 \in S-I$ and hence $(a_1,a_3) \subseteq S-I$.

To show the reverse containment first note that $S-I \subseteq S$ since $0 \in I$.  Now, let $x \in S\backslash (a_1,a_3)$.  Then we know that $x = t_2a_2 + t_4a_4$ for some $t_2,t_4 \ge 0$.  Note that if $t_2 \ge q_3$, then \\

\begin{center}
\begin{tabular}{ll}

$x$\hspace{-.09in} &= $(t_2 - q_3)a_2 + q_3a_2 + t_4a_4$\\
    &= $(t_2 - q_3)a_2 + q_2a_3 + t_4a_4$ \ \ (by 1.7 (5))\\ 
    
\end{tabular}
\end{center}

\vspace{.1in}

Thus $x \in (a_1,a_3)$.  Similarly, if $t_4 \ge q_1$, then \\

\begin{center}
\begin{tabular}{ll}

$x$\hspace{-.09in} &= $t_2a_2 + (t_4-q_1)a_4 + q_1a_4$\\
    &= $t_2a_2 + (t_4-q_1)a_4 + q_4a_1$ \ \ (by 1.7 (5))\\ 
    
\end{tabular}
\end{center}

\vspace{.1in}

Thus $x \in (a_1,a_3)$.  Therefore, we may assume that $0 \le t_2 \le q_3-1$ and $0 \le t_4 \le q_1-1$.  We now examine a very specific quantity.\\

\noindent
\begin{tabular}{ll}

$(q_3-1)a_2 + (q_1-1)a_4 + n$\hspace{-.09in} &= $(q_3-1)a_2 + (q_1-1)a_4 + a_2-a_1$ \vspace{.04in}\\
 & $= q_3a_2 + (q_1-1)a_4 - a_1$ \vspace{.04in}\\
 & $= q_3a_2 + (q_1-1)a_4 + (q_1-2)a_1 - (q_1-1)a_1$ \vspace{.04in}\\
 & $= g(T) - (q_1-1)a_1$ \ \ \ (by (2.2)) \vspace{.04in}\\
 & $= g(S)$ \ \ \ (by (2.6)) \vspace{.04in}\\
                               
\end{tabular}

\vspace{.2in}

This means $(q_3-1)a_2 + (q_1-1)a_4 + n \notin S$.  Since $S$ is closed under addition, we conclude $t_2a_2 + t_4a_4 + n \notin S$ for $0 \le t_2 \le q_3-1$ and $0 \le t_4 \le q_1-1$.  It follows that $x+n \notin S$ and thus $x \notin S-I$.  We now know $S-I \subseteq (a_1,a_3)$ whence the two relative ideals are equal.\\

We have now completed the proof of (2.1).\\

\noindent
{\bf (2.9) Example:}\ Based upon the investigation in this section, the following pairs of semigroups and relative ideals are perfect $2 \times 2$ bricks:\\

\begin{center}

$(\langle 24,25,35,36\rangle,(0,1))$\\

 \vspace{.1in}
$(\langle 15,22,33,40\rangle,(0,7))$\\

\vspace{.1in}
$(\langle 28,45,81,98\rangle,(0,17))$\\

\end{center}

\vspace{.2in}

To establish the existence of infinitely many perfect $2 \times 2$ bricks it suffices to prove that there are infinitely many unitary numerical semigroups.  Consider $S = \langle 2(2z+1),5z,5(z+1),3(2z+1)\rangle$.  It is quick to check that if $z \ge 3$ and $5 \not|\ 2z+1$, then $S$ is a unitary numerical semigroup.  There are, in fact, infinitely many examples of these ``families'' of unitary numerical semigroups.\\

\vspace{.2in}
%%%%%%%%%%%%%%%%%%   Section 3   %%%%%%%%%%%%%%%%%%

\begin{center}

{\bf 3. Other Examples, Notes and Open Questions}\\

\end{center}

This investigation into unitary numerical semigroups and $k \times m$ bricks began as a search for bricks beyond the three known examples given at the end of \cite{HS}.  The search employed a ``brute force'' computer program written in $C++$.  We searched numerical semigroups of the form $S=\langle s_1,\dots ,s_t\rangle$ where $2 \le t \le 5$ and $2 \le s_i \le 50$ for each $i$.  The relative ideals were of the form $I = (u_1,\dots ,u_l)$ where $2 \le \mu_S(I) \le 1+\frac{t}{2}$ and $0 \le u_j \le g(S)-m(S)$ for each $j$.  The smallest generator for $I$ was always $0$.  The search took approximately 36 hours on a 3 GHz Pentium 4 machine running Linux.  The search revealed over 12,000 examples of numerical semigroups and relative ideals that form bricks of various dimensions.  A small sample of these is listed below.  \\

\noindent
{\bf (3.1) Example:}\\

\begin{tabular}{ll}

$S = \langle 10,15,18,27\rangle$ &  \\
$I=(0,2)$ \ \  $S-I=(18,25)$ & $2 \times 2$ brick \vspace{.1in}\\
$S = \langle 21,28,36,48\rangle$ &  \\
$I=(0,13)$ \ \ $S-I=(36,56,63)$ & $2 \times 3$ brick \vspace{.1in}\\
$S = \langle 15,17,21,24,27\rangle$ &  \\
$I=(0,8)$ \ \ $S-I=(24,30,34,36)$ & $2 \times 4$ brick \vspace{.1in}\\
$S = \langle 21,24,38,39\rangle$ &  \\
$I=(0,4,6)$ \ \  $S-I=(72,77,80)$ & $3 \times 3$ brick \vspace{.1in}\\
$S = \langle 27,30,36,44\rangle$ &  \\
$I=(0,1,6)$ \ \ $S-I=(87,98,101,110)$ & $3 \times 4$ brick\\

\end{tabular}

\ \\

We offer some observations about the bricks that were found during this search.\\

\noindent
{\bf (3.2) Observations:}\\

\noindent
(1)  We found no perfect bricks in dimensions other than $2 \times 2$.\\

\noindent
(2)  We found no bricks with multiplicity 9, 11 or 13.\\

\noindent
(3)  Let $(S,I)$ be any $2 \times 2$ brick.  Recall this means $\mu_S(I)=\mu_S(S-I)=2$ and $\mu_S(I+(S-I))=4$.  Let $I = (0,n)$ and $S-I = (a_1,a_3)$.  If we define $\hat{S}$ to be the numerical semigroup minimally generated by $\{a_1,a_2=a_1+n,a_3,a_4=a_3+n\}$ and let $\hat{I}$ be the relative ideal of $\hat{S}$ minimally generated by $\{0,n\}$, then $\hat{S}-\hat{I}=(a_1,a_3)$ and $(\hat{S},\hat{I})$ turns out to be a perfect $2 \times 2$ brick.\\

We finish this preliminary investigation into unitary numerical semigroups and $k \times m$ bricks by offering the following open questions.\\

\noindent
{\bf (3.3) Open Questions:}\\

\noindent
(1) By (2.1) we know that every unitary numerical semigroup yields a perfect $2 \times 2$ brick.  Is the converse true?  That is, if $(S,I)$ is a perfect $2 \times 2$ brick, is it true that $S$ must be unitary?\\

\noindent
(2) Is it true that every imperfect $2 \times 2$ brick will yield a perfect $2 \times 2$ brick in the manner described in (3.2 (3))?\\

\noindent
(3) Do there exist any perfect bricks that are not $2 \times 2$?\\

\noindent
(4) Do there exist bricks of all possible dimensions?  That is, given $k > 1$ and $m > 1$, does there exist a $k \times m$ brick?\\

\vspace{.1in}

{}     


\begin{thebibliography}{9}
\bibitem{A} M. Auslander, {\sl Modules over unramified regular local rings}, Illinois J. Math {\bf 5} (1961), 631-647

\bibitem{BDF} V. Barucci, D. Dobbs and M. Fontana, {\sl Maximality properties in numerical semigroups and applications to one-dimensional analytically irreducible local domains}, Memoirs of the American Mathematical Society {\bf 125} (1997).

\bibitem{C} P. Constapel, {\sl Vanishing of Tor and torsion in tensor products}, Comm. Algebra {\bf 24} (1996), 833 - 846

\bibitem{FGH} R. Froberg, C. Gottlieb and R. Haggkvist, {\sl On numerical semigroups}, Semigroup Forum {\bf 35} (1987), 63 - 83

\bibitem{H1} K. Herzinger, {\sl Torsion in the tensor product of an ideal with its inverse}, Comm. Algebra {\bf 24} (1996), 3065-3083

\bibitem{H3} K. Herzinger, {\sl The number of generators for an ideal and its
dual in a numerical semigroup}, Comm. Algebra {\bf 27} (1999), 4673-4689

\bibitem{HS} K. Herzinger and R. Sanford, {\sl Minimal generating sets for relative ideals in numerical semigroups of multiplicity eight}, Comm. Algebra {\bf 32} (2004), 4713-4731

\bibitem{HW1}  C. Huneke and R. Wiegand, {\sl Tensor products of modules and the rigidity of Tor}, Math. Ann. {\bf 299} (1994), 449-476

\bibitem{HW2}  C. Huneke and R. Wiegand, {\sl Tensor products of modules, rigidity and local cohomology}, Math. Scand. {\bf 81} (1997), 161-183

\bibitem{J} S. M. Johnson, {\sl A linear Diophantine problem}, Canad. J. Math. {\bf 12} (1960), 390-398

\bibitem{K} E. Kunz, {\sl The value-semigroup of a one-dimensional Gorenstein ring}, Proc. Amer. Math. Soc. {\bf 25} (1970), 748-751

\bibitem{MCH} M. Madero-Craven and K. Herzinger, {\sl Apery sets of numerical semigroups}, Comm. Algebra (to appear)

\bibitem{S} J. J. Sylvester, {\sl Mathematical questions with their solutions}, Educational Times {\bf 41} (1884), 21 

\end{thebibliography}
\end{document}